\documentclass[12pt]{article}
\usepackage{amssymb, amsfonts, amsmath, amscd, amsthm}
\usepackage[all]{xy}

\numberwithin{equation}{section}

\newcommand{\tD}{\operatorname{\widetilde{\mathcal{D}}_q}}
\newcommand{\tM}{\operatorname{\widetilde{\mathcal{M}}}}

\newcommand{\tU}{{\widetilde{U}_q}}
\newcommand{\tTX}{\operatorname{\widetilde{T}^\star X}}
\newcommand{\Ind}{\operatorname{Ind}}

\newcommand{\Hom}{\operatorname{Hom}}
\newcommand{\End}{\operatorname{End}}

\newcommand{\mood}{\hbox{\ensuremath{
\operatorname{mod}}}}

\newcommand{\ad}{\hbox{\ensuremath{\operatorname{ad}}}}
\newcommand{\mO}{\mathcal{O}}

\newcommand{\Z}{\ensuremath{\mathcal{Z}}}
\newcommand{\Zl}{\ensuremath{\mathcal{Z}^{(l)}}}
\newcommand{\ZHC}{\ensuremath{\mathcal{Z}^{HC}}}

\newcommand{\C}{\ensuremath{{\mathbb{C}}}}
\newcommand{\D}{\ensuremath{\mathcal{D}_{G_q}}}
\newcommand{\Dl}{\ensuremath{\mathcal{D}^\lambda_q}}

\newcommand{\g}{\ensuremath{\mathfrak{g}}}

\newcommand{\n}{\ensuremath{\mathfrak{n}}}
\newcommand{\h}{\ensuremath{\mathfrak{h}}}

\newcommand{\bb}{\ensuremath{\mathfrak{b}}}

\newcommand{\Oq}{\ensuremath{\mathcal{O}_q}}

\newcommand{\UA}{\ensuremath{U_\mathcal{A}}}
\newcommand{\UAres}{\ensuremath{U^{res}_\mathcal{A}}}

\newcommand{\Fp}{\ensuremath{\mathbb{F}_p}}
\newcommand{\UBr}{\ensuremath{U^{res}_q(\bb)}}
\newcommand{\W}{\ensuremath{\mathcal{W}}}
\newcommand{\Tl}{\ensuremath{T^\star X^{\lambda}}}
\newcommand{\M}{\ensuremath{\mathcal{M}}}
\newcommand{\MBG}{\ensuremath{\mathcal{M}_{B_q}(G_q)}}

\newcommand{\DBG}{\ensuremath{\mathcal{D}^\lambda_{B_q}(G_q)}}
\newcommand{\tDBG}{\ensuremath{\widetilde{\mathcal{D}}_{B_q}(G_q)}}

\newcommand{\U}{\ensuremath{U^{fin}_q}}
\newcommand{\Ul}{\ensuremath{U^{\lambda}_q}}

\newcommand{\A}{\ensuremath{\mathcal{A}}}

\newcommand{\OA}{\ensuremath{\mathcal{O}_\mathcal{A}}}

\theoremstyle{plain}

\newtheorem*{Thm*}{Theorem}
\newtheorem*{Cor*}{Corollary}
\newtheorem{Thm}{Theorem}[section]
\newtheorem{Prop}[Thm]{Proposition}
\newtheorem{Lem}[Thm]{Lemma}

\newtheorem{Cor}[Thm]{Corollary}

\theoremstyle{definition}
\newtheorem{defi}[Thm]{Definition}

\theoremstyle{remark}
\newtheorem{Rem}[Thm]{Remark}

\begin{document}
\title{Proof of the De Concini-Kac-Procesi conjecture}

\author{Kobi Kremnizer}
\maketitle

\begin{abstract}
In this paper we prove a conjecture by De Concini, Kac and Procesi \cite{CP} (Corollary \ref{conj}):

\begin{Cor*}
 The dimension of any $M\in U_q-\mood^\chi$ is divisible by $l^{codim_\mathcal{B}\mathcal{B}_\chi}$.
\end{Cor*}

\end{abstract}

\section{Introduction}
In this paper we prove a conjecture by De Concini, Kac and Procesi \cite{CP} (Corollary \ref{conj}):

\begin{Cor*}
 The dimension of any $M\in U_q-\mood^\chi$ is divisible by $l^{codim_\mathcal{B}\mathcal{B}_\chi}$.
\end{Cor*}
The conjecture follows from the following theorem (Theorem \ref{transpoly}):

\begin{Thm*}
Fix $\chi$ a nilpotent element in the big cell and a regular $\lambda\in P$. For any module $M\in U_q-\mood^{\lambda,\chi}$ there exists a polynomial $\mathbf{d}_M\in \frac{1}{R}\mathbb{Z}[P^*]$
of degree $\leqslant dim(\mathcal{B}_\chi)$, such that for any $\mu\in P$ in the closure of the alcove of $\lambda$, we have
\begin{equation}
dim(T_\lambda^\mu(M))=\mathbf{d}_M(\mu)
\end{equation}
 $\mathbf{d}_M(\mu)=l^{dim(\mathcal{B})}\mathbf{d}_M^0(\frac{\mu+\rho}{l})$ for another polynomial $\mathbf{d}_M^0\in \frac{1}{R}\mathbb{Z}[P^*]$, such that $\mathbf{d}_M^0(\mu)\in\mathbb{Z}$ for $\mu\in P$.
 \end{Thm*}
 
The proof of this theorem follows quite easily from the results of \cite{BK2} and some usual arguments concering the Euler characteristic.

It must be emphasized that \cite{BMR} prove a similar result in positive characteristic (the Kac-Weisfeiler conjecture first proved by Premet) using geometric methods. This allowed us to adapt their proof to the quantum case by using the geometric construction of modules for the quantum group from \cite{BK2}. Thus this paper in another indication of the close relationship between positive characteristic representations and root of unity representations. This relationship will be investigated more in the future. 

We would like to thank Dmitriy Rumynin that suggested we should prove this conjecture.

\thanks{The author was supported in part by NSF grant DMS-0602007.}
 
\section{Quantum groups} 
\subsection{Conventions}\label{Conventions}

Let $\C$ be the field of complex numbers and fix $q \in \C^\star$.

\textbf{We always assume that if $q$ is a root of unity it is
primitive of odd order and in case $G$ has a component of type
$G_2$ the order is also prime to 3.} Let $\A$ be the local ring
$\mathbb{Z}[\nu]_{\mathfrak{m}}$, where $\mathfrak{m}$ is the
maximal ideal in $\mathbb{Z}[\nu]$ generated by $\nu - 1$ and a
fixed odd prime $p$.

\subsection{Root data} Let $\g$ be a semi-simple Lie
algebra and let $\h \subset \bb$ be a Cartan subalgebra contained
in a Borel subalgebra of $\g$. Let $R$ be the root system, $\Delta
\subset R_+ \subset R$ a basis and the positive roots. Let $P
\subset \h^\star$ be the weight lattice and $P_+$ the positive
weights; the $i$'th fundamental weight is denoted by $\omega_i$
and $\rho$ denotes the half sum of the positive roots. Let $Q
\subset P$ be the root lattice and $Q_+ \subset Q$ those elements
which have non-negative coefficients with respect to the basis of
simple roots. Let $\mathcal{W}$ be the Weyl group of $\g$. We let
$<\,,\,>$ denote a $\mathcal{W}$-invariant bilinear form on
$\h^\star$ normalized by $<\gamma,\gamma> = 2$ for each short root
$\gamma$.

Let $T_P = \Hom_{groups}(P,k^\star)$ be the character group of $P$
with values in $k$ (we use additive notation for this group). If
$\mu \in P$, then $<\mu,P> \subset \mathbb{Z}$ and hence we can
define $q^\mu \in T_P$ by the formula $q^\mu(\gamma) =
q^{<\mu,\gamma>}$, for $\gamma \in P$. If $\mu \in P, \lambda \in
T_P$ we write $\mu + \lambda = q^\mu + \lambda$. Note that the
Weyl group naturally acts on $T_P$.

\subsection{Quantized enveloping algebra $U_q$ and quantized
algebra of functions $\Oq$.} Let $U_q$ be the simply connected
quantized enveloping algebra of $\g$ over $\C$. Recall that $U_q$
has algebra generators $E_\alpha, F_\alpha, K_\mu$, $\alpha,
\beta$ are simple roots, $\mu \in P$ subject to the relations
\begin{equation}\label{qgr1}
K_{\lambda} K_{\mu} = K_{\lambda+\mu},\;\; K_0 = 1,
\end{equation}
\begin{equation}\label{qgr2}K_\mu E_\alpha K_{-\mu} = q^{<\mu, \alpha>} E_\alpha,\;\; K_\mu F_\alpha K_{-\mu} = q^{-<\mu, \alpha>}
F_\alpha,\end{equation}
\begin{equation}\label{qgr3}[E_\alpha,\, F_\beta] = \delta_{\alpha, \beta} {{K_\alpha - K_{-\alpha}}\over{q_\alpha-q^{-1}_\alpha}}\end{equation}
and certain Serre-relations that we do not recall here. Here

$q_\alpha = q^{d_\alpha}$, $d_\alpha = <\alpha,\alpha>/2$. (We
have assumed that $q^2_\alpha \neq 1$.)

Let $G$ be the simply connected algebraic group with Lie algebra
$\g$, $B$ be a Borel subgroup of $G$ and $N \subset B$ its
unipotent radical. Let $\bb = \operatorname{Lie} B$ and
$\mathfrak{n} = \operatorname{Lie} N$ and denote by $U_q(\bb)$ and
$U_q(\mathfrak{n})$ the corresponding subalgebras of $U_q$. Then
$U_q(\bb)$ is a Hopf algebra, while $U_q(\mathfrak{n})$ is only an
algebra. Let $\Oq = \Oq(G)$ be the algebra of matrix coefficients
of finite dimensional type-1 representations of $U_q$. There is a
natural pairing $(\;,\;): U_q \otimes \Oq \to \C$. This gives a
$U_q$-bimodule structure on $\Oq$ as follows
\begin{equation}\label{23}
ua = a_1 (u,a_2),\;\; au = (u,a_1)a_2, \,\, u \in U_q, a \in \Oq
\end{equation}
Then $\Oq$ is the (restricted) dual of $U_q$ with respect to this
pairing. We let $\Oq(B)$ and $\Oq(N)$ be the quotient algebras of
$\Oq$ corresponding to the subalgebras $U_q(\bb)$ and
$U_q(\mathfrak{n})$ of $U_q$, respectively, by means of this
duality. Then $\Oq(B)$ is a Hopf algebra and $\Oq(N)$ is only an
algebra.

There is a braid group action on $U_q$. For each $w \in \W$, we
get an automorphism $T_w$ of $U_q$.

\subsection{Integral versions of $U_q$.} Let
$\UAres$ be the Lusztig's integral form of $U_q$, the $\A$-algebra
in $U_q$ generated by divided powers $E^{(n)}_\alpha =
E^n_\alpha/{[n]_{d_\alpha}!}$, $F^{(n)}_\alpha =
F^n_\alpha/{[n]_{d_\alpha}!}$, $\alpha$ a simple root, $n \geq 1$
(where $[m]_d = {\prod^m_{s=1} {q^{d\cdot s} -  q^{-d\cdot
s}}\over{q^{d} -  q^{-d}}}$) and the $K_\mu$'s, $\mu \in P$. There
is also the De Consini-Kac integral form $\UA$, which is generated
over $\A$ by the $E_\alpha, F_\alpha$ and $K_\mu$'s. The
subalgebra $\UA$ is preserved by the adjoint action of $\UAres$:
$ad_{\UAres} (\UA) \subset \UA$. The operators $T_w$ from section
2.1.2 preserves the integral versions.

$\OA$ is defined to be the dual of $\UAres$. This is an $\A$-sub
Hopf algebra of $\Oq$.

\subsection{Finite part of $U_q$.} The algebra $U_q$ acts on
itself by the adjoint action $\ad: U_q \to U_q$ where $\ad(u)(v) =
u_1 v S(u_2)$. Let $U^{fin}_q$ be the finite part of $U_q$ with
respect to this action:
$$ U^{fin}_q = \{v \in U_q; \dim \ad(U_q)(v) < \infty\}.$$
This is a subalgebra. (See \cite{JL}.)

We can also give an integral version of the finite part as the
finite part of the action of $\UAres$ on $\UA$.Thus by
specializing we get a subalgebra of $U_q$ for every $q$. Of
course, when specialized to generic $q$ this coincides with the
previous definition.

\subsection{Specializations and Frobenius maps.}
For any ring map $\phi: \A \to R$ we put $U_R = \UA \otimes_\A R$
and $U^{res}_R = \UAres \otimes_\A R$. If $R = \C$ and $\phi(\nu)
= q$, there are three different cases: $q$ is a root of unity, $q
= 1$ and $q$ is generic. Then $U_R = U_q$.

There is the also the ring map $\A \to \Fp$, sending $\nu \to 1$.
Then $U_{\Fp} /(K-1) = U(\g_p)$, the enveloping algebra of the Lie
algebra $\g_p$ in characteristic $p$.

For any $\UA$-module (resp. $\UAres$-module) $M_\A$ we put $M_R =
M_\A \otimes_\A R$. This is an $U_R$-module (resp.
$U^{res}_R$-module). When $R = \C$ we simply write $M = M_\C$.

When $q$ is a root of unity, we have the Frobenius map: $U^{res}_q
\to U(\g)$. Its algebra kernel is denoted by $\mathfrak{u}_q$. We
also have the Frobenius map $U^{res}_q(\bb) \to U(\bb)$, with
algebra kernel $b_q$. These maps induces dual maps $\mO = \mO(G)
\hookrightarrow \Oq$ and $\mO(B) \hookrightarrow \Oq(B)$.

For each $q$ there exists a map $U_q \to U^{res}_q$ whose image is
$\mathfrak{u}_q$ and whose algebra kernel is $\Zl$ (see section
2.1.7 below for the definition of $\Zl$).

\subsection{Verma modules.}\label{Verma} For each $\lambda \in T_P$ there is the
one dimensional $U_q(\bb)$-module $\C_\lambda$ which is given by
extending $\lambda$ to act by zero on the $E_\alpha$'s. The
Verma-module $M_\lambda$ is the $U_q$-module induced from
$\C_\lambda$. If $\mu \in P$ we write $M_\mu = M_{q^\mu}$. A point
important for us is that $M_\lambda$ carries an
$U^{res}_q(\bb)$-module structure defined as follows:
$U^{res}_q(\bb)$ acts on $U_q$ by restricting the adjoint action
of $U^{res}_q$ on $U_q$. This induces a $U_q(\bb)$-action on the
quotient $M_\lambda$ of $U_q$. Since this action is locally finite
it corresponds to an $\Oq(B)$-comodule action on $M_\lambda$.
\textbf{Note!!} As a $U^{res}_q(\bb)$-module $M_\lambda$ has
\textit{trivial highest weight}. (In case $q$ is generic,
$U^{res}_q(\bb) = U_q(\bb)$ and then the $U^{res}_q(\bb)$ action
on $M_\lambda$ described above is the same as the $U(\bb)$-action
on $M_\lambda \otimes \C_{-\lambda}$.)

Verma module $M_\lambda$ has an integral version $M_{\lambda,\A}$.

\subsection{Centers of $U_q$ and definition of $\tU$.} Let $\Z$ denote the center of $U_q$.
When $q$ is a $p$'th root of unity $\Z$ contains the
Harish-Chandra center $\ZHC$ and the $l$'th center $\Zl$ which is
generated by the $E^{l}_\alpha$, $F^{l}_\alpha$, $K^{l}_\mu$ and
$K^{-l}_\mu$'s. In fact, $\Z = \Zl \otimes_{\Zl \cap \ZHC} \ZHC$.
There is the Harish-Chandra homomorphism $\ZHC \to \mO(T_P)$ that
maps isomorphically to the $W$-invariant even part of
$\mO(T_{P})$. We define $\tU = U_q \otimes_{\ZHC} \mO(T_P)$.

\subsection{Some conventions.} We shall frequently refer to a right
(resp. left) $\Oq$-comodule as a left (resp. right) $G_q$-module,
etc. If we have two right $\Oq$-comodules $V$ and $W$, then
$V\otimes W$ carries the structure of a right $\Oq$-comodule via
the formula
$$
\delta(v\otimes w) = v_1\otimes w_1 \otimes v_2w_2$$ We shall
refer to this action as the {\it tensor} or {\it diagonal} action.
A similar formula exist for left comodules.

\section{Quantum flag variety}
Here we recall the definition and basic properties of the quantum
flag variety from \cite{BK1}.

\subsubsection{Category $\MBG$.} The composition
\begin{equation}\label{e2}
\Oq \to \Oq \otimes \Oq \to \Oq \otimes \Oq(B)\end{equation}
defines a right $\Oq(B)$-comodule structure on  $\Oq$. A
$B_q$-equivariant sheaves on $G_q$ is a triple $(F, \alpha,
\beta)$ where $F$ is a vector space, $\alpha: \Oq \otimes F \to F$
a left $\Oq$-module action and $\beta: F \to F \otimes \Oq(B)$ a
right $\Oq(B)$-comodule action such that $\alpha$ is a right
comodule map, where we consider the tensor comodule structure on
$\Oq(G) \otimes F$.
\begin{defi}\label{hejsan} We denote $\MBG$ to be the category of $B_q$-equivariant
sheaves on $G_q$. Morphisms in $\MBG$ are those compatible with
all structures.\end{defi}

If $q = 1$, the category $\M_B(G)$ is equivalent to the category
$\M(G/B)$ of quasi-coherent sheaves on $\mathcal{B}=G/B$.

\begin{defi} We define the induction functor $\Ind: B_q-\mood$ to $\MBG$, $\Ind
V = \Oq \otimes V$ with the tensor $B_q$-action and the
$\Oq$-action on the first factor. For $\lambda \in P$ we define a
line bundle $\Oq(\lambda) = \Ind \C_{-\lambda}$.
\end{defi}
\begin{defi} The global section functor $\Gamma: \MBG \to \C-\mood$
is defined by
$$\Gamma(M) = \Hom_{\MBG}(\Oq, M) = \{m \in M; \Delta_B(m) =
m\otimes 1\}.$$ This is the set of $B_q$-invariants in $M$.
\end{defi}

The category $\MBG$ has enough injectives, so derived functors
$R\Gamma$ are well-defined. We showed that $R^i\Gamma(\Ind V) =
H^i(G_q/B_q,V)$, where $H^i(G_q/B_q,\;)$ is the $i$'th derived
functor of the functor $V \to \Gamma(\Ind V)$ from $B_q-\mood$ to
$\C-\mood$.

We proved a quantum version of Serre's basic theorem on projective
schemes: Each $M \in \MBG$ is a quotient of a direct sum of
$\Oq(\lambda)$'s and each surjection $M \twoheadrightarrow M'$ of
noetherian objects in $\MBG$ induces a surjection
$\Gamma(M(\lambda)) \twoheadrightarrow \Gamma(M'(\lambda))$ for
$\lambda >> 0$.

Here the notation $\lambda >> 0$ means that $<\lambda,
\alpha^\wedge>$ is a sufficiently large integer for each simple
root $\alpha$ and $M(\lambda) = M \otimes \C_{-\lambda}$ is the
$\lambda$-twist of $M$.

Let $V \in G_q-\mood$. Denote by $V \vert B_q$ the restriction of
$V$ to $B_q$ and by $V^{triv}$ the trivial $B_q$-module whose
underlying space is $V$. We showed that $\Ind V \vert B_q$ and
$\Ind V^{triv}$ are isomorphic in $\MBG$. In particular
\begin{equation}\label{Gmod}
\Gamma(\Ind V \vert B_q) = V \vert B_q \otimes \Gamma(\Oq) = V
\vert B_q, \hbox{ for } V \in G_q-\mood
\end{equation}
\subsubsection{ \MBG at a root of unity}
In case $q$ is a root of unity we have the following Frobenius
morphism:
\begin{equation}
Fr_*:\MBG \rightarrow \M(G/B)
\end{equation}
defined as
\begin{equation}
 N \mapsto N^{b_q}
\end{equation}

Using the description of $\MBG$ as $Proj(A_q)$ where $A_q=
\bigoplus V_{q,\lambda}$ \cite{BK1} and similarly $\M(G/B)=
Proj(A)$ where $A= \bigoplus V_\lambda$, we see that $Fr_*$ is
induced from the quantum Frobenius map $A\hookrightarrow A_q$. It
follows that:
\begin{Prop}\label{Frob} $Fr_*$ is exact and faithful.
\end{Prop}

 This functor has a left adjoint which we denote by $Fr^*$.

 We will need the following:
 \begin{Lem}
  $Fr^*(\mathcal{F}(\lambda))\simeq Fr^*(\mathcal{F})(l\lambda)$
 \end{Lem}
\begin{Lem}
 $Fr_*(Fr^*(\mathcal{F})(\lambda))\simeq \mathcal{F}\otimes Fr_*(\mathcal{O}_{q,\lambda})$
\end{Lem}

 \begin{Lem}\label{wellknown} $Fr_*(\mathcal{O}_q(-\rho))\simeq \mathcal{O}(-\rho)^{l^{dim(\mathcal{B})}}$
 \end{Lem}

  \section{The ring $\D$, the category of $\Dl$-modules and translation functors} 
 \subsection{$\Dl$-modules}
  We need the following important
\begin{Rem}\label{Integral versions.} All objects described in the preceding chapters are
defined over $\A$. For any specialization $\A \to R$ and any
object $Obj$ we denote by $Obj_R$ its $R$-form. For the functors
we don't use any subscripts; so, for instance, there is the
functor $\Ind: B_{q,R}-\operatorname{mod} \to \MBG_R$.
\end{Rem}
Recall the $U_q$-bimodule structure on $\Oq$ given by \ref{23}.
Now, as we have two versions of the quantum group we pick the
following definition of the ring of differential operators on the
group (the $crystalline$ version).
\begin{defi} We define the ring of quantum differential operators on $G_q$ to be
the smash product algebra $\D := \Oq \star U_q$. So $\D = \Oq
\otimes U_q$ as a vector space and multiplication is given by
\begin{equation}\label{e71}
a \otimes u \cdot b \otimes v = au_1(b)\otimes u_2v.
\end{equation}
\end{defi}
We consider now the ring $\D$ as a left $U^{res}_q$-module, via
the left $U^{res}_q$-action on $\Oq$ in \ref{23} and the left
adjoint action of $U^{res}_q$ on $U_q$; this way $\D$ becomes a
module algebra for $U^{res}_q$: In the following we will use the
restriction of this action to $U^{res}_q(\bb) \subset U^{res}_q$.
As $U_q$ is not locally finite with respect to the adjoint action,
this $U^{res}_q(\bb)$-action doesn't integrate to a $B_q$-action.
Thus $\D$ is not an object of $\MBG$; however, $\D$ has a
subalgebra $\D^{fin} = \Oq \star \U$ which belongs to $\MBG$. This
fact will be used below.

\begin{defi}\label{d12} Let $\lambda \in T_P$. A $(B_{q},\lambda)$-equivariant $\D$-module
is a triple $(M, \alpha, \beta)$, where $M$ is a $\C$-module,
$\alpha: \D \otimes M \to M$ a left $\D$-action and $\beta^{res}:
M \to M\otimes \Oq(B)$ a right $\Oq(B)$-coaction. The latter
action induces an $U^{res}_q(\bb)$-action on $M$ again denoted by
$\beta^{res}$. We have the natural map $U_q(\bb) \to \UBr$ which
together with $\beta^{res}$ gives an action $\beta$ of $U_q(\bb)$
on $M$. We require
\smallskip

\noindent $i)$ The $U_q(\bb)$-actions on $M \otimes \C_\lambda$
given by $\beta \otimes \lambda$ and by $(\alpha\vert_{U_q(\bb)})
\otimes \operatorname{Id}$ coincide.

\noindent $ii)$ The map $\alpha$ is $U^{res}_q(\bb)$-linear with
respect to the $\beta$-action  on $M$ and the action on $\D$.

\smallskip

These objects form a category denoted $\DBG$. There is the
forgetful functor $\DBG \to \MBG$. Morphisms in $\DBG$ are
morphisms in $\MBG$ that are $\D$-linear.
\end{defi}
We defined $\Dl$ as the maximal quotient of $\D$ which is an
object of $\DBG$ and showed that
\begin{equation}\label{relevantstruct1}
\Dl = \Ind M_\lambda
\end{equation}
as an object in $\DBG$. (See section \ref{Verma} for the
$B_q$-action $=$ $U^{res}_q(\bb)$-action on $M_\lambda$.). The
global section functor $\Gamma: \DBG \to \M$ is the functor of
taking $B_q$ invariants (with respect to the action $\beta$); we
have $\Gamma = \Hom_{\DBG}(\Dl,\;)$.

Hence, in particular $\Gamma(\Dl) = \End_{\DBG}(\Dl)$ (which
explains the ring structure on $\Gamma(\Dl)$).

Let $q$ be a $l$-th root of unity. In \cite{BK2} the following were proved:
\begin{Prop}\label{RGammaD} We have $i)$ $\Ul \cong R\Gamma(\Dl)$
and $ii)$ $\tU \cong R\Gamma(\tD)$ (if $l$ is a prime $p >$
Coxeter number of $G$.)
\end{Prop}
\begin{Prop}\label{Azumaya}
 $\Dl$ is an
Azumaya algebra  over a dense subset of $\Tl$.
\end{Prop}
 \begin{Prop}\label{Azsplit} The action map $U^{fin}_q \otimes_Z O_{\tTX^{twist}} \to
\tD $ induces an isomorphism $U^{fin}_q \otimes_Z
O_{\tTX^{twist,unram}} \cong \tD\mid _{\tTX^{twist,unram}} $.
\end{Prop}
\begin{Thm}\label{Daffthm}
$R\Gamma: D^b(\DBG) \to D^b(\Ul-\mood)$ is an equivalence of
categories.
\end{Thm}

\subsection{Translation Functors} 
For $\lambda,\mu\in P$ define the translation functor $T_\lambda^\mu:U_q^\lambda-\mood\to U_q^\mu-\mood$
$T_\lambda^\mu(M)=[V_{\mu-\lambda}\otimes M]_\mu$. Here $V_{\mu-\lambda}$ is the standard $G_q$-module with an extremal weight $\mu-\lambda$. 
This functor is also defined on the larger categories $U_q-\mood^\lambda$ of modules with a generalized central character.

 We want to define translation functors on $\mathcal{D}$-modules as well. We first define a category $\tDBG$ of
$\D$-modules that contains all $\Dl$, $\lambda \in T_P$, a
"torsor".
\begin{defi}\label{d1t} An object of $\tDBG$
is a triple $(M, \alpha, \beta)$, where $\alpha: \D \otimes M \to
M$ a left $\D$-action and $\beta^{res}: M \to M\otimes \Oq(B)$ a
right $\Oq(B)$-coaction.
\smallskip

\noindent $i)$ The $U_q(\n_+)$-actions on $M$ given by $\beta\vert
U_q(\n_+)$ and by $(\alpha\vert_{U_q(\n_+)})$ coincide.

\noindent $ii)$ The map $\alpha$ is $U_q(\bb)$-linear with respect
to the $\beta$-action  on $M$ and  the action on $\D$.
\end{defi}

The functor of $B_q$-invariants (global sections) will be denoted $\widetilde{\Gamma}$.
Put $\tM = U_q(\g)/\sum_{\alpha \in R_+} U_q(\g) \cdot E_\alpha$
(a "universal" Verma module) and define
\begin{equation}\label{tildeD}
\tD = \Oq \otimes \tM
\end{equation}
$\tD$ inherits an $U^{res}_q(\bb)$-module structure from $\D$, so
$\tD$ is an object in $\tDBG$.

We define a category $\tDBG_\lambda$ by also requiring that a power of the ideal defined by $\lambda$
kills the module.
We will denote the functor of $B_q$-invariants (global sections) by $\widetilde{\Gamma}_\lambda$.

Note that we have a pair of exact functors:
\begin{equation}
 i_\lambda:\tDBG_\lambda\rightleftarrows\tDBG:[ ]_\lambda
\end{equation}
where $i_\lambda$ is the natural inclusion and $[]_\lambda$ is the projection.

For any $B_q$-module $V$ and any $\DBG$-module $M$ we can define a $\tDBG$-module $M\otimes V$
where we just twist the $B_q$ action and let $\D$ act on the first factor.
We have the following equalities:
\begin{equation}
T_\lambda^\mu(\widetilde{R\Gamma}_\lambda(M)=[\widetilde{R\Gamma}_\lambda M\otimes V_{\mu-\lambda}]_\mu=[\widetilde{R\Gamma}(M\otimes V_{\mu-\lambda})]_\mu\simeq \widetilde{R\Gamma}_\mu([M\otimes V_{\mu-\lambda}]_\mu) 
\end{equation}

Note also that that $M\otimes V_{\mu-\lambda}$ is an extension of terms $(M\otimes k_{-\nu})\otimes V_{\mu-\lambda}(\nu)$ (see \cite{BK1} last section) and hence $[M\otimes V_{\mu-\lambda}]_\mu$
is an extension of terms $(M\otimes k_{-\nu})\otimes V_{\mu-\lambda}(\nu)$ with $\nu\in \mathcal{W}_\lambda^\mu-\lambda\subseteq\mathcal{W}_{\lambda-\mu}$.

If $\lambda,\mu$ lie in the same closed alcove then
\begin{equation}
\mathcal{W}_\lambda^\mu=(\lambda+\mathcal{W}_{\mu-\lambda})\cap W_{aff}\bullet\mu=(W_{aff})_\lambda\bullet\mu\subseteq \lambda+W\cdot(\mu-\lambda) 
\end{equation}

and if $\mu$ is also in the closure of the facet of $\lambda$ then since $(W_{aff})_\mu\subseteq(W_{aff})_\lambda$ we get that $\mathcal{W}_\lambda^\mu=\lbrace\mu\rbrace$.
Hence we get from the previous discussion:
\begin{Lem}\label{down}
 Under the above conditions on $\mu,\lambda$ we have
 \begin{equation}
  T_\lambda^\mu(\widetilde{R\Gamma}_\lambda M)\simeq \widetilde{R\Gamma}_\mu(M\otimes k_{\lambda-\mu})
 \end{equation}

\end{Lem}

\section{Dimensions of $U_{q,\chi}$-modules}
In this section we follow \cite{BMR}. Set $R=\prod_\alpha<\rho,\check{\alpha}>$ where $\alpha$ runs over the set of positive roots.
\begin{Thm}\label{transpoly} Fix $\chi$ a nilpotent element in the big cell and a regular $\lambda\in P$. For any module $M\in U_q-\mood^{\lambda,\chi}$ there exists a polynomial $\mathbf{d}_M\in \frac{1}{R}\mathbb{Z}[P^*]$
of degree $\leqslant dim(\mathcal{B}_\chi)$, such that for any $\mu\in P$ in the closure of the alcove of $\lambda$, we have
\begin{equation}
dim(T_\lambda^\mu(M))=\mathbf{d}_M(\mu)
\end{equation}
 $\mathbf{d}_M(\mu)=l^{dim(\mathcal{B})}\mathbf{d}_M^0(\frac{\mu+\rho}{l})$ for another polynomial $\mathbf{d}_M^0\in \frac{1}{R}\mathbb{Z}[P^*]$, such that $\mathbf{d}_M^0(\mu)\in\mathbb{Z}$ for $\mu\in P$.
 \end{Thm}
 \begin{Rem}
  Without loss of generality we can assume that $\chi\in\mathfrak{n}$ (Using the infinite dimensional group of automorphisms of \cite{CKP}). 
 \end{Rem}

 \begin{Cor}\label{conj}
 The dimension of any $M\in U_q-\mood^\chi$ is divisible by $l^{codim_\mathcal{B}\mathcal{B}_\chi}$
 \end{Cor}
 \textit{proof}
 We can choose a regular $\lambda$ and $\mu$ in the closure of the $\lambda$-facet such that $N=T_\lambda^\mu(M)$. Then by theorem \ref{transpoly} $dim(N)=l^{dim\mathcal{B}}\mathbf{d}_M^0(\frac{\mu+\rho}{l})$. But then since $R$ is prime to $p$ for
 $p>h$ we get that the denominator of $\mathbf{d}_M^0(\frac{\mu+\rho}{l})$ is a power of $p$ less than
 $deg(\mathbf{d}^0_M)\leq dim(\mathcal{B}_\chi)$ and hence $dim(N)$ is divisible by a power of $p$ greater than $dim(\mathcal{B})-deg(\mathbf{d}^0_M)\geq codim(\mathcal{B}_\chi)$.
\begin{Lem}\label{splitbun}
Let $\mathcal{M}_\chi^\lambda$ be the splitting bundle for the restriction of the Azumaya algebra $\Dl$ to $\mathcal{B}_\chi$. We have an equality in $K^0(\mathcal{B}_\chi)$:
\begin{equation}
[\mathcal{M}_\chi^\lambda]=[Fr_*\mathcal{O}_q(l\rho+\lambda)\vert_{\mathcal{B}_\chi}]
\end{equation}
\end{Lem}   
\textit{proof}
As $\mathcal{M}_\chi^\lambda=\mathcal{M}_\chi^0(\lambda)$ by the construction of the splitting it is enough to check for $\lambda=-\rho$.
On the zero section (that is on $\mathcal{B}_0$) the Azumaya algebra is split with splitting bundle
$\mathcal{O}_\mathcal{B}^{l^{dim\mathcal{B}}}\simeq Fr_*(\mathcal{O}_q((l-1)\rho))$. 
We also have that $U_q^{-\rho}\vert_\mathfrak{n}\simeq End(M_{-\rho})$, hence
if we denote by $\widetilde{\mathfrak{n}}$ the preimage of $\mathfrak{n}$ under the Springer map then there exists a splitting bundle $\widetilde{\mathcal{M}}$ whose restriction to $\mathcal{B}_\chi$ is $\mathcal{M}_\chi^{-\rho}$. Now we can use the map $\mathcal{B}_\chi\times\mathbb{A}^1\to \widetilde{\mathfrak{n}}$ given by $(x,t)\mapsto (x,t\chi)$ to pull back this splitting bundle $\widetilde{\mathcal{M}}$. By the rational invariance of $K^0$ we get the result. $\Box$ 
\begin{Lem}\label{polynomial}
For any $\mathcal{F}\in D^b(Coh(\mathcal{B}))$ or $\mathcal{F}\in D^b(Coh(\mathcal{B}_q))$ there exists a polynomial $\mathbf{d}_\mathcal{F}\in \frac{1}{R}\mathbb{Z}[P^*]$ such that for $\lambda\in P$ the Euler characteristic of $R\Gamma(\mathcal{F}\otimes\mathcal{O}(\lambda))$
equals $\mathbf{d}(\lambda)$. We have in the case of $\mathcal{B}$
\begin{equation}\label{dim}
deg(\mathbf{d}_\mathcal{F})\leqslant dim( supp(\mathcal{F})).
\end{equation}
If $\mathcal{F}\in D^b(Coh(\mathcal{B}))$
\begin{equation}\label{identity}
\mathbf{d}_{Fr^*(\mathcal{F})}(\mu)=l^{dim\mathcal{B}}\mathbf{d}_\mathcal{F}(\frac{\mu+(1-l)\rho}{l})
\end{equation}
\end{Lem}
\textit{proof} The first statement is true for line bundles by the Weyl character formula (and its quantum version) and the fact that the $K$-group is generated by classes of line bundles (also in the quantum case). The inequality follows from Grothendieck-Riemann-Roch (\cite{BMR}).
To prove \ref{identity} it is enough to check it for $\mu=l\nu-\rho$. Then
\begin{equation}
 Fr_*(Fr^*(\mathcal{F})(l\nu-\rho))\simeq Fr_*(Fr^*(\mathcal{F}(\nu))(-\rho))\simeq\mathcal{F}(\nu)\otimes Fr_*(\mathcal{O}(-\rho))
\end{equation}
but from \ref{wellknown} we get
\begin{equation}
 \mathcal{F}(\nu)\otimes Fr_*(\mathcal{O}(-\rho))\simeq\mathcal{F}(\nu-\rho)^{\oplus l^{dim(\mathcal{B})}}.
\end{equation}
$\Box$

\subsection{proof of the theorem}
Let $\mathcal{F}_M\in Coh(\widetilde{\mathcal{B}}_\chi)$ be such that $\mathcal{L}^\lambda M\simeq \mathcal{F}_M\otimes\mathcal{M}_\lambda$. The existence of such follows from the localization theorem and the Azumaya splitting over te fibers. Let $[\mathcal{F}_M]\in K(Coh(\widetilde{\mathcal{B}}_\chi))=K(\mathcal{B}_\chi)$ be its class. By Lemma \ref{down} we have
\begin{equation}
T_\lambda^\mu(M)=R\Gamma(\mathcal{L}^\lambda(M)(\mu-\lambda))=R\Gamma(\mathcal{F}_M\otimes\mathcal{M}_\lambda(\mu-\lambda))=R\Gamma(\mathcal{F}_M\otimes\mathcal{M}_\mu)
\end{equation}
Denote by $\int$ the Euler charcteristic of $R\Gamma$.  We have
\begin{equation}
dim(T_\lambda^\mu(M))=\int_{\mathcal{B}_\chi}[\mathcal{F}_M]\cdot[\mathcal{M}_\mu]
\end{equation}
where the multiplication is the action of of $K^0$ on $K$.
By Lemma \ref{splitbun} we have 
\begin{equation}
\int_{\mathcal{B}_\chi}[\mathcal{F}_M]\cdot i^*[Fr_*\mathcal{O}_{q,l\rho+\mu}]=
\int_\mathcal{B}i_*[\mathcal{F}_M]\cdot[Fr_*\mathcal{O}_{q,l\rho+\mu}]=
\int_{\mathcal{B}_q}Fr^*(i_*[\mathcal{F}_M(l\rho+\mu)])
\end{equation}
where $i:\mathcal{B}_\chi\hookrightarrow\mathcal{B}$ is the inclusion.
By Lemma \ref{polynomial} we get
\begin{equation}
 dim(T_\lambda^\mu)=\mathbf{d}_{Fr^*(i_*\mathcal{F}_M)}(l\rho+\mu)=l^{dim\mathcal{B}}\mathbf{d}_{i_*\mathcal{F}_M}(\frac{\mu+\rho}{l}).
\end{equation}
So by \ref{dim} and \ref{identity} the polynomial $\mathbf{d}_{i_*\mathcal{F}_M}$ has thr required properties.


\begin{thebibliography}{}
\bibitem[ABG]{ABG} S. Arkhipov, R. Bezrukavnikov, V. Ginzburg,
\textit{Quantum groups, the loop Grassmannian, and the Springer resolution.},
J. Amer. Math. Soc.  17  (2004), 595-678.









\bibitem[BK1]{BK1} E. Backelin and K. Kremnizer
\textit{Quantum flag varieties, equivariant quantum
$\mathcal{D}$-modules, and localization of Quantum groups.}, Adv. in Math.
203 (2006) 408-429.

\bibitem[BK2]{BK2} E. Backelin and K. Kremnizer
\textit{Localization for quantum groups at a root of unity.}, math.RT/0407048.





\bibitem[BMR]{BMR} A. Bezrukavnikov, I. Mirkovic and D. Rumynin,
\textit{Localization for a semi-simple Lie algebra in prime
characteristic.}, arXiv:math.RT/0205144.





\bibitem[CKP]{CKP} C. De Concini, V.G. Kac and C. Procesi,
\textit{Quantum coadjoint action }, JAMS v.5, number 1, (1992)
151-189.

\bibitem[CP]{CP} C. De Concini, C. Procesi,
\textit{Quantum groups}, LNM 1565, 31-140.




\bibitem[JL]{JL} A. Joseph, G. Letzter
\textit{Local finiteness for the adjoint action for quantized
enveloping algebras}, J. Algebra 153, (1992), 289-318.
















\end{thebibliography}
 \end{document}